\documentclass[letterpaper, 10 pt, conference]{ieeeconf}  

\IEEEoverridecommandlockouts                              

\title{\LARGE \bf
Pareto Control Barrier Function for Inner Safe Set Maximization Under Input Constraints
}

\author{Xiaoyang Cao$^{1}$, Zhe Fu$^{2,*}$ and Alexandre M. Bayen$^{2}$
\thanks{*Corresponding author. Email: zhefu@berkeley.edu}
\thanks{$^{1}$Tsinghua University, Beijing, China}%
\thanks{$^{2}$University of California, Berkeley, CA, USA}%
}

\usepackage{amsthm}
\usepackage[T1]{fontenc}
\usepackage{cite}
\usepackage{amssymb}
\usepackage{amsmath}
\usepackage{algorithm}
\usepackage{algpseudocode}
\usepackage{booktabs}
\usepackage[hidelinks]{hyperref}

\newtheoremstyle{bolditalicdefinition} 
  {\topsep}   
  {\topsep}   
  {\itshape}  
  {}          
  {\bfseries} 
  {.}         
  {.5em}      
  {}          

\theoremstyle{bolditalicdefinition}
\newtheorem{definition}{Definition} 

\newtheoremstyle{bolditalicproblem} 
  {\topsep}   
  {\topsep}   
  {\itshape}  
  {}          
  {\bfseries} 
  {.}         
  {.5em}      
  {}          

\theoremstyle{bolditalicproblem}
\newtheorem{problem}{Problem} 

\newtheorem{lemma}{Lemma}

\usepackage{graphicx}

\begin{document}

\maketitle
\thispagestyle{empty}
\pagestyle{empty}

\begin{abstract}
Control Barrier Functions (CBFs) enforce safety in dynamical systems by ensuring trajectories stay within prescribed safe sets. However, traditional CBFs often overlook realistic input constraints, potentially limiting their practical effectiveness. To address this research gap, this paper introduces the Pareto Control Barrier Function (PCBF) algorithm, which employs a Pareto multi-task learning framework to simultaneously balance the competing objectives of maintaining safety and maximizing the volume of the safe set under input limitations. The PCBF algorithm is computationally efficient and scalable to high-dimensional systems. The effectiveness of PCBF is demonstrated through comparisons with Hamilton-Jacobi reachability on an inverted pendulum and extensive simulations on a 12-dimensional quadrotor system, where PCBF consistently outperforms existing methods by yielding larger safe sets while ensuring robust safety under input constraints. Codes are available at \url{https://github.com/XiaoyangCao1113/Pareto_CBF}.
\end{abstract}

\section{INTRODUCTION}
Control Barrier Functions (CBFs) provide a powerful method to ensure the safety of autonomous systems by maintaining system trajectories within a predefined safe set \cite{ames2019control}. They have been successfully applied in various domains, including robotics \cite{nguyen2015safety} and autonomous vehicles \cite{alan2023control}.

However, ensuring safety becomes challenging when systems operate under input constraints \cite{agrawal2021safe, breeden2021high, liu2023safe, so2024train}. Traditional CBFs typically assume the availability of arbitrary control inputs \cite{lindemann2018control, xu2018safe}, which allows the system to be guided within the safe set without limitations. Yet, this assumption often ignores the fact that control inputs are constrained by actuator limits, energy resources, and other system-specific factors. As a result, the safe set may no longer be forward invariant under such restricted inputs \cite{agrawal2021safe}.

Recent research has tried to address this problem. Agrawal and Panagou \cite{agrawal2021safe} introduced Input Constrained CBFs (ICCBFs), iteratively defining functions to form an inner safe set \( C^* \) that is forward invariant under input constraints. Breeden and Panagou \cite{breeden2021high} introduced Zeroing CBFs (ZCBFs) for high-relative-degree systems, ensuring forward invariance by converting higher-order constraints into ZCBFs. However, both approaches tend to be overly conservative, as they do not guarantee the largest possible safe set, potentially limiting the system's operational flexibility.

Neural CBFs (NCBFs) have emerged as an effective approach to ensure safety in complex, high-dimensional systems, where traditional CBFs may be difficult to design \cite{dawson2023safe}. These methods have been extended to handle parametric uncertainties, unknown dynamic environments, and safe navigation tasks \cite{richards2018lyapunov, kolter2019learning, xiao2023barriernet}. However, few works have integrated input constraints into NCBFs. Liu et al. proposed the non-saturating CBF (NSCBF) to handle input constraints solving a min-max optimization \cite{liu2023safe}. Though effective, the complexity of this approach leads to longer training times. Oswin et al. introduced a policy NCBF based on policy iteration \cite{so2024train}. However, their method does not guarantee the maximal inner safe set either.

Hamilton-Jacobi (HJ) reachability analysis theoretically computes the largest control-invariant set \cite{mitchell2005time}. Recent work has combined HJ reachability with CBFs to mitigate the overly conservative problem of HJ policies \cite{choi2021robust}. However, HJ reachability is computationally intensive and practically limited to systems with fewer than five dimensions, restricting its applicability to higher-dimensional dynamical systems.

To address these limitations, we make the following contributions:
\begin{itemize}
    \item We propose the PCBF algorithm to determine the largest inner safe set under input constraints.
    \item We enhance training efficiency by applying a Gaussian sampling method to the inner safe set's interior.
    \item We validate the proposed approach through comparisons with HJ reachability on an inverted pendulum and demonstrate its scalability via simulations on a 12-dimensional quadrotor system.
\end{itemize}

The remainder of this article is organized as follows: Section \ref{sec:preliminaries} provides preliminaries, including problem definition, CBFs, Neural CBFs, and Pareto Multi-task Learning. Section \ref{sec:methodology} presents our methodology, detailing the PCBF algorithm. Section \ref{sec:experiments} demonstrates the effectiveness of our approach through experiments on an inverted pendulum and a quadrotor system. Section \ref{sec:conclusions} concludes the paper and discusses future work directions.

\section{PRELIMINARIES}

In this section, we present the theoretical framework of CBFs and Neural CBFs, and how they ensure safety guarantees.

\label{sec:preliminaries}
\subsection{Problem Definition}
Consider continuous-time, control-affine dynamics:

\begin{equation}
    \dot{x} = f(x) + g(x)u,
    \label{eq:system_dynamics}
\end{equation}

where $x \in \mathcal{X} \subseteq \mathbb{R}^{n}$, $u \in \mathcal{U} \subseteq \mathbb{R}^{m}$, and $f, g$ are locally Lipschitz continuous functions. To ensure safety, the concept of a safe set is introduced.

\begin{definition}[Safe Set]
The \textit{Safe Set} \( S \) is defined as:
\begin{equation}
    S \triangleq \{x \in \mathcal{X} : h(x) \geq 0\},
\end{equation}
where \( h: \mathcal{X} \to \mathbb{R} \) is a continuously differentiable function.
\end{definition}

The safe set \( S \) is defined to include all states $x$ in which the system remains safe during operation. In practice, the safe set often comes from the physical constraints and safety requirements inherent in the system, such as avoiding collisions or maintaining stability. However, the presence of input constraints makes it challenging to ensure forward invariance of the entire safe set \cite{agrawal2021safe}. Thus, the concept of an inner safe set is introduced.

\begin{definition}[Forward Invariant]
A set \( C \subseteq \mathbb{R}^n \) is \textit{forward invariant} under a dynamical system \eqref{eq:system_dynamics} if for any initial condition \( x(0) \in C \), the trajectory \( x(t) \in C \) for all \( t \geq 0 \).
\end{definition}

\begin{definition}[Inner Safe Set \cite{agrawal2021safe}]
An \textit{Inner Safe Set} \( C \) is a non-empty, closed subset of the safe set \( S \) such that there exists a feedback controller \( \pi: C \to \mathcal{U} \) that renders \( C \) forward invariant under the system dynamics.
\end{definition}

This research aims to find the largest such inner safe set given the system's input constraints. The problem is formulated as below:

\begin{problem}[Inner Safe Set Maximization]
Given a dynamical system \eqref{eq:system_dynamics} and a safe set \( S \subseteq \mathcal{X} \), determine the largest possible inner safe set \( C^* \subseteq S \) and a control policy \( \pi: C^* \to \mathcal{U} \) such that \( C^* \) remains forward invariant under the system dynamics.
\label{problem:inner-safe-set-max}
\end{problem}

\subsection{Control Barrier Functions}

Control Barrier Functions (CBFs) provide a systematic approach to ensuring the forward invariance of safe sets.

\begin{definition}[Control Barrier Functions (CBFs)]
    A function \( h: \mathcal{X} \to \mathbb{R} \) is considered a CBF if there exists an extended class-\(\mathcal{K}\) function     \footnote{An extended class-\(\mathcal{K}\) function is a strictly increasing, continuous function \( \alpha: \mathbb{R} \to \mathbb{R} \) with \( \alpha(0) = 0 \).}  \( \alpha \)
    such that for all states \( x \in \mathcal{X} \), the following inequality holds:
    \begin{equation}
    \sup_{u \in \mathcal{U}} \left[ L_f h(x) + L_g h(x) u \right] \geq -\alpha(h(x)),
    \label{def-cbf}
    \end{equation}
    where \( L_f h(x) = \nabla h(x)^\top f(x) \) and \( L_g h(x) = \nabla h(x)^\top g(x) \) are the Lie derivatives of \( h(x) \) along the vector fields \( f(x) \) and \( g(x) \), respectively.
\end{definition}

To compute a safe control input $u(x)$, the following CBF Quadratic Program (CBF-QP) \cite{garg2024learning} is formulated.
The objective is to find a control input that remains as close as possible to the nominal control policy $\pi_{nom}(x): \mathcal{X} \rightarrow \mathcal{U}$ while enforcing the safety constraints. According to Theorem 2 in \cite{ames2019control}, the solution to problem \ref{eq:cbf-qp} guarantees the forward invariance of the safe set.

\begin{problem}[CBF-QP]
\begin{equation}
    \begin{aligned}
        \min_{u \in \mathcal{U}} & \quad \| u - \pi_{nom}(x) \|^2 \\
        \text{s.t.} & \quad L_f h(x) + L_g h(x) u \geq -\alpha(h(x)).
    \end{aligned}
\end{equation}
\label{eq:cbf-qp}
\end{problem}

The nominal policy $\pi_{nom}$ is derived using unconstrained control design methods (e.g., LQR, MPC, etc.) that do not inherently guarantee safety.
For instance, consider a self-driving car following a planned route. 
Its nominal policy is the ideal steering and acceleration commands for smooth lane tracking. 
However, if a pedestrian suddenly steps onto the road, the CBF-QP adjusts the commands from $\pi_{nom}$ just enough to safely avoid a collision, while still trying to stay close to the original plan.

When \( \mathcal{U} \) is unbounded, 
the affine structure of \eqref{def-cbf} in \( u \) ensures that a valid CBF and a corresponding control input \( u \) can always be found, provided \( L_g h(x) \not\equiv 0 \).
However, under input constraints, the control input \( u \) that satisfies \eqref{def-cbf} may no longer remain within the admissible set \( \mathcal{U} \). As a result, certain states within the safe set may lose forward invariance. This raises the key challenge of identifying the largest inner safe set that is compatible with input constraints.

\subsection{Neural CBFs}

Neural CBFs (NCBFs) address the difficulty of hand-designing CBFs for complex systems by using neural networks to approximate a CBF and its corresponding controller \cite{dawson2023safe, garg2024learning}.
However, such NCBFs still require a predefined safe set. 
If this safe set contains any unsafe states, the neural networks will fail to converge.
Liu et al. \cite{liu2023safe} proposed the NSCBF $h_{\theta}(x)$, which removes the need for a predefined safe set and incorporates input constraints into its formulation.

\begin{equation}
    h_{\theta}(x) = h(x) - \left( \text{nn}_{\theta}(x) - \text{nn}_{\theta}(x_e) \right)^2,
    \label{eq:def-ncbf}
\end{equation}

\noindent where \( \text{nn}_{\theta}(x) \) is a neural network with parameter $\theta \in \mathbb{R}^p$, and \( x_e \) is a known safe state.
Typically, in problems aimed at maintaining the system's proximity to a stable point, \( x_e \) itself is the stable point under some appropriate control law. For obstacle avoidance problems, \( x_e \) can be a point far from the obstacle.
This method ensures that the learned inner safe set \( C_{h_{\theta}} = \{ x \in \mathcal{X}: h_{\theta}(x) \geq 0 \} \) is always a subset of the original safe set. Also, \( C_{h_{\theta}} \) is always a non-empty set, guaranteeing that there are sufficient sample points during the training process.

Although this method has shown success in specific cases, it has several limitations. First, it requires sampling on \( \partial C_{h_{\theta}} \) to calculate the saturation risk, which is computationally expensive; second, it does not guarantee the identification of the largest possible safe set. Moreover, even for the maximal inner safe set, the loss function remains non-zero, which can lead to suboptimal solutions.

\section{METHODOLOGY}
\label{sec:methodology}
The Pareto multi-task learning framework \cite{zitzler1999multiobjective} is introduced into traditional NCBFs to develop the PCBF algorithm. This approach maintains a balance between ensuring safety and maximizing the inner safe set under input constraints. 
The NCBF formulation from equation \eqref{eq:def-ncbf} is adopted throughout this study.

\subsection{Loss Function Design}

Two loss functions have been proposed for effective training of the PCBF: the feasibility loss and the volume loss. Computing these two losses requires sampling in the state space. A multivariate Gaussian sampling method is employed to ensure sufficient data density near \(x_e \in \mathcal{X}\): $\mathbb{X}_{\mathcal{N}} = \{x_i \mid x_i \overset{\text{i.i.d.}}{\sim} \mathcal{N}(x_e, \Sigma), \ i = 1, 2, \dots, N_{\mathcal{N}}\} \subset \mathcal{X}$, where \( N_{\mathcal{N}} \) is the number of sampled points. The covariance matrix \( \Sigma \) is diagonal, with each diagonal element defined as \( \sigma_{ii} = \frac{\Delta x_i}{k} \), where \( \Delta x_i \) is the state space length in the \( i \)-th dimension, and \( k \) is a positive hyperparameter.

\textbf{Feasibility loss}: The feasibility loss \( \mathcal{L}_{\text{feas}} \) is necessary for the forward invariance of the learned inner safe set \( C_{h_{\theta}} \):

\begin{equation}
    \mathcal{L}_{\text{feas}} = \frac{1}{N_{\theta}} \sum_{x_i \in \mathbb{X}_{\mathcal{N}} \cap C_{h_{\theta}}} \max \{ 0, -\sup_{u \in \mathcal{U}} \phi_{\theta}(x_i, u) \},
\end{equation}

\begin{equation}
    \phi_{\theta}(x, u) = L_f h_{\theta}(x) + L_g h_{\theta}(x) u + \alpha(h_{\theta}(x)),
\end{equation}

\noindent where $N_{\theta}$ is the number of sampling points in $\mathbb{X}_{\mathcal{N}} \cap C_{h_{\theta}}$. For simplicity, this study assumes linear input constraints, making \( \mathcal{U} \) a polyhedron. Since \( \phi_{\theta}(x, u) \) is affine in \( u \), the supremum over \( \mathcal{U} \) is necessarily achieved at some vertex of \( \mathcal{U} \), allowing closed-form computation of $\mathcal{L}_{\text{feas}}$.

\begin{equation}
    \sup_{u \in \mathcal{U}}\phi_{\theta}(x, u) = \max_{u \in \mathcal{V}(\mathcal{U})} \phi_{\theta}(x, u),
\end{equation}

\noindent where $\mathcal{V}(\mathcal{U})$ is the set of vertexes of $\mathcal{U}$.

\textbf{Volume loss}: The volume loss \( \mathcal{L}_{\text{vol}} \) is designed to encourage the expansion of the learned inner safe set by minimizing the deviation from $\mathcal{S}$:

\begin{equation}
    \mathcal{L}_{\text{vol}} = \frac{1}{N_{\theta}} \sum_{x_i \in \mathbb{X}_{\mathcal{N}} \cap C_{h_{\theta}}} (\text{nn}_{\theta}(x_i) - \text{nn}_{\theta}(x_e))^2,
\end{equation}

However, a simple linear combination of these loss functions does not guarantee a maximized inner safe set. The feasibility loss and the volume loss are inherently competing as both cannot be zero simultaneously: Expanding the safe set violates safety constraints, and satisfying constraints reduces the set volume.
To address this challenge, the Pareto multi-task learning framework is introduced to balance these conflicting objectives.

\begin{figure}[ht]
    \centering
    \includegraphics[width=.7\linewidth]{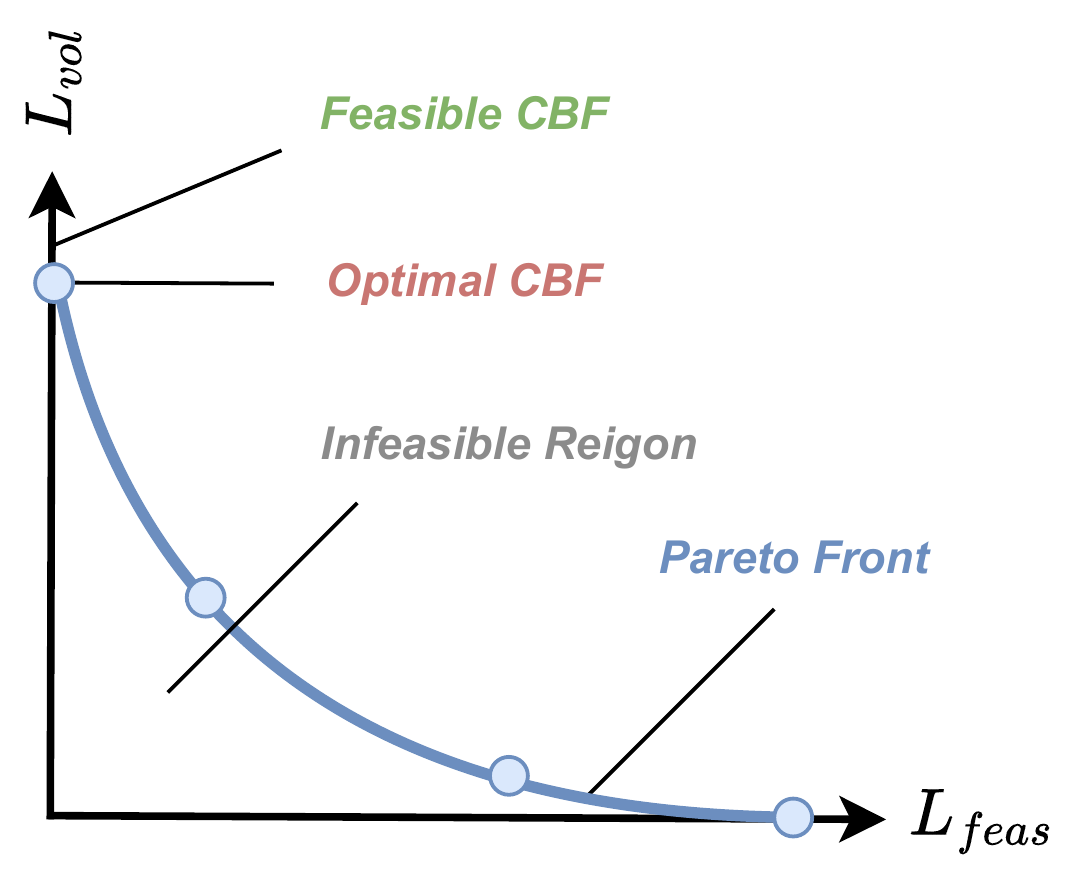}
    \caption{Pareto front of the Inner Safe Set Maximization Problem. The blue line represents the Pareto front, which delineates the boundary between the feasible and infeasible regions in the objective space. No solution $\theta$ can produce objective values $[\mathcal{L}_{\text{feas}}(\theta) \quad \mathcal{L}_{\text{vol}}(\theta)]^T$ that lie below this curve. The vertical axis represents solutions with zero feasibility loss. All feasible CBFs lie on this axis. The optimal CBF, with parameter $\theta^*$, lies at the intersection of the Pareto front and the vertical axis.}
    \label{fig:sketch-pareto}
\end{figure}

\subsection{Pareto Multi-task Learning (PMTL)}
Consider \( M \) correlated tasks with a loss vector:

\begin{equation}
    \mathbf{L}(\theta) = \begin{bmatrix}
        \mathcal{L}_1(\theta) & \mathcal{L}_2(\theta) & \cdots & \mathcal{L}_M(\theta)
    \end{bmatrix}^T, \theta \in \mathbb{R}^p,
\end{equation}

\noindent where \( \mathcal{L}_i(\theta) \) is the loss of the \( i \)-th task. In general, it is impossible to find a single solution \( \theta \) that simultaneously optimizes all the objectives.
Instead, a collection of Pareto optimal solutions is identified, each representing a distinct optimal trade-off among the various objectives \cite{zitzler1999multiobjective}.

\begin{definition}[Pareto Dominance \cite{lin2019pareto}]
    A solution \( \theta \) dominates another solution \( \hat{\theta} \) if:
    \begin{itemize}
    \item $\forall i \in \{1, \dots, M\}, \, \mathcal{L}_i(\theta) \leq \mathcal{L}_i(\hat{\theta}),$
    \item $\exists j \in \{1, \dots, M\} \text{ such that } \mathcal{L}_j(\theta) < \mathcal{L}_j(\hat{\theta}).$
    \end{itemize}
\end{definition}

\begin{definition}[Pareto Optimality \cite{lin2019pareto}]
A solution \( \theta \) is Pareto optimal if no other solution \( \hat{\theta} \) exists that dominates \( \theta \).
\end{definition}

In other words, a solution \( \theta \) is considered Pareto optimal if no other feasible solution \( \hat{\theta} \) exists that can reduce any loss without simultaneously increasing at least one other loss.

\begin{definition}[Pareto Front \cite{lin2019pareto}]
The Pareto Front \( P_{\mathcal{L}} \) is the set of all loss vectors associated with the Pareto optimal solutions:
\begin{equation}
    P_{\mathcal{L}} = \{\mathbf{L}(\theta) \mid \theta \in \mathbb{R}^p \text{ is Pareto optimal}\}.
\end{equation}
\end{definition}

The Pareto front represents the set of optimal trade-offs among competing objectives achievable by any solution. PMTL provides a framework for simultaneously optimizing multiple objectives and converges to the Pareto front. Two common PMTL methods are:

\textbf{Linear combination}: A non-negative linear weighted sum is utilized to aggregate the losses of all tasks into a single loss: $\mathcal{L}(\theta) = \sum_{i=1}^M \lambda_i \mathcal{L}_{i}(\theta)$, where $\lambda_i \geq 0$ is the weight for the $i$-th task. Although straightforward, this approach has some limitations: first, selecting appropriate values for $\lambda_i$ is challenging and often requires extensive empirical tuning; second, this method can only produce solutions on the convex portion of the Pareto front \cite{boyd2004convex}.

\textbf{Gradient-based method}: Fliege and Svaiter \cite{fliege2000steepest} proposed a gradient-based method generalizing the single objective steepest descent algorithm. The algorithm's update rule is given by: $\theta_{t+1} = \theta_t + \eta d_t$, where $\eta > 0$ is the learning rate and the search direction $d_t \in \mathbb{R}^p$ is obtained by solving following optimization problems:

\begin{figure}[ht]
    \centering
    \includegraphics[width=0.75\linewidth]{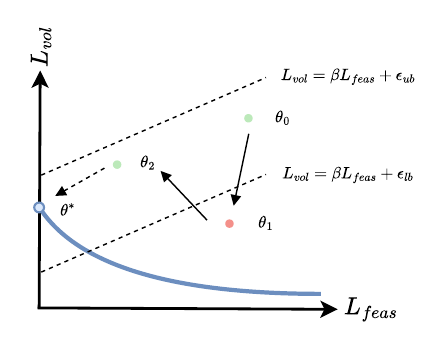}
    \caption{The blue solid curve represents the Pareto front, while the black dashed lines indicate the upper and lower bounds. The green point denotes a feasible solution (i.e., a solution $\theta$ with $\mathbf{L}(\theta)$ lies between the two bounds), the red point represents an infeasible solution, and the blue point corresponds to $\theta^*$, the solution of problem \ref{eq:problem-ism-loss}. If the update causes $\theta$ to fall outside the feasible region, the PCBF algorithm will gradually guide $\theta$ back into the feasible region bounded by the upper and lower bounds.}
    \label{fig:pcbf}
\end{figure}

\begin{equation}
\begin{aligned}
(d_t, \alpha_t) = \arg & \min_{d \in \mathbb{R}^p, \alpha \in \mathbb{R}} \alpha + \frac{1}{2}\|d\|^2 \\
\text{s.t.} & \quad \nabla \mathcal{L}_i(\theta_t)^T d \leq \alpha, \quad i = 1, \ldots, M.
\end{aligned}
\label{eq:problem-base}
\end{equation}

The solutions of this problem satisfy the following lemma:

\begin{lemma}[Fliege and Svaiter \cite{fliege2000steepest}]
Let $(d_t, \alpha_t)$ be the solution of problem \eqref{eq:problem-base}. Then:
\begin{enumerate}
    \item If $\theta_t$ is Pareto critical, then $d_t = 0 \in \mathbb{R}^p$ and $\alpha_t = 0$.
    \item If $\theta_t$ is not Pareto critical, then
    \begin{align}
        \alpha_t &\leq -\frac{1}{2}\|d_t\|^2 < 0, \label{eq:alpha-bound} \\
        \nabla \mathcal{L}_i(\theta_t)^T d_t &\leq \alpha_t, \quad i = 1, \ldots, M.
    \end{align}
\end{enumerate}
\label{lemma:steepest-descent}
\end{lemma}

A solution $\theta$ is said to be Pareto critical if no solution in its neighborhood can perform better in all tasks. If $\theta_t$ is Pareto critical, $d_t=0$ and $\theta_t$ will not be updated. If $\theta_t$ is not Pareto critical, then $d_t$ is a descent direction for all tasks.

\textbf{Connection to problem \ref{problem:inner-safe-set-max}}: Problem \ref{problem:inner-safe-set-max} can be reformulated as a PMTL problem with $M=2$ and 
$\mathbf{L}(\theta) = [\mathcal{L}_{\mathrm{feas}}(\theta) \quad \mathcal{L}_{\mathrm{vol}}(\theta)]^T$ emperically.

\begin{problem}[Reformulation of Inner Safe Set Maximization]
\begin{equation}
\begin{aligned}
    \theta^* = \arg & \min_{\theta \in \mathbb{R}^p} \quad \mathcal{L}_{\mathrm{vol}}(\theta) \\
    \text{s.t.} & \quad \mathcal{L}_{\mathrm{feas}}(\theta) = 0.
\end{aligned}
\end{equation}
\label{eq:problem-ism-loss}
\end{problem}

The solution of problem \ref{eq:problem-ism-loss}, denoted as \(\theta^*\), achieves zero feasibility loss while minimizing the volume loss, and lemma \ref{lemma:ism-pareto} shows that $\theta^*$ is on the Pareto front. Figure \ref{fig:sketch-pareto} provides a visual explanation of the trade-offs between the two losses on the Pareto front. 

\begin{lemma}
Let $\theta^*$ be the solution to problem \ref{eq:problem-ism-loss}, then \(\theta^*\) is Pareto optimal, and the loss vector \([\mathcal{L}_{\mathrm{feas}}(\theta^*) \mathcal{L}_{\mathrm{vol}}(\theta^*)]^T\) lies on the Pareto front.
\label{lemma:ism-pareto}
\end{lemma}

The PCBF algorithm is proposed to effectively solve problem \ref{eq:problem-ism-loss} by applying gradient-based PMTL methods. 

\subsection{Pareto Control Barrier Function (PCBF)}
Initially, we attempted to employ a linear combination of $\mathcal{L}_{\text{feas}}$ and $\mathcal{L}_{\text{vol}}$ as the total loss function. 
However, this approach often led to one of two outcomes:
either \(\mathcal{L}_{\text{vol}} \approx 0\) (resulting in a model that fails to remove potentially hazardous states) 
or \(\mathcal{L}_{\text{feas}} \approx 0\) but with a large \(\mathcal{L}_{\text{vol}}\) (yielding an excessively small inner safe set).

To mitigate these issues, the objective space is partitioned into three regions using two parallel lines (as shown in Figure \ref{fig:pcbf}). 
The intermediate region is designed to contain \(\mathbf{L}(\theta^*)\). 
The motivation behind the PCBF algorithm is to restrict \(\mathbf{L}(\theta)\) within this intermediate region, thereby avoiding convergence to the aforementioned extreme outcomes while guiding the gradient descent towards \(\theta^*\).
Specifically, we define an upper bound:
$
\mathcal{L}_{\text{vol}}(\theta) = \beta \mathcal{L}_{\text{feas}}(\theta) + \epsilon_{\text{ub}},
$
and a lower bound:
$
\mathcal{L}_{\text{vol}}(\theta) = \beta \mathcal{L}_{\text{feas}}(\theta) + \epsilon_{\text{lb}},
$
where \(\beta\), \(\epsilon_{\text{ub}}\), and \(\epsilon_{\text{lb}}\) are hyperparameters.
To ensure $\mathbf{L}(\theta^*)$ lies between these bounds, constraints are imposed
\footnote{Although the exact value of $\mathcal{L}_{\text{vol}}(\theta^*)$ is unknown, we can obtain an estimate by first performing a pre-training process using the linear combination loss.}
such that $\beta > 0$ and $0 \leq \epsilon_{\text{lb}} \leq \mathcal{L}_{\text{vol}}(\theta^*) \leq \epsilon_{\text{ub}}$.

In iteration \(t+1\), the parameter update is given by:
$\theta_{t+1} = \theta_t + \eta d_t$,
where \(d_t\) is the update direction and $\eta >0$ is the learning rate. To ensure that the updated parameters remain within the specified region, \(d_t\) is determined based on the position of \(\mathbf{L}(\theta_t)\) in the objective space:

\begin{enumerate}
    \item If $\mathcal{L}_{\text{vol}}(\theta_t) > \beta \mathcal{L}_{\text{feas}}(\theta_t) + \epsilon_{\text{ub}}$:
    \begin{equation}
    \begin{aligned}
        (d_t, \alpha_t) = \arg & \min_{d \in \mathbb{R}^p, \alpha \in \mathbb{R}} \; \alpha + \frac{1}{2} \|d\|_2^2 \\
        \text{s.t.} & \quad [\nabla (\mathcal{L}_{\text{vol}}(\theta_t) - \beta \mathcal{L}_{\text{feas}}(\theta_t))]^T d \leq \alpha, \\
        & \quad \nabla \mathcal{L}_{\text{feas}}(\theta_t)^T d \leq \alpha.
    \end{aligned}
    \label{eq:primal_1}
    \end{equation}

    \item If $\beta \mathcal{L}_{\text{feas}}(\theta_t) + \epsilon_{\text{lb}} \leq \mathcal{L}_{\text{vol}}(\theta_t) \leq \beta \mathcal{L}_{\text{feas}}(\theta_t) + \epsilon_{\text{ub}}$:
    \begin{equation}
    \begin{aligned}
        \label{eq:primal_2}
        (d_t, \alpha_t) = \arg & \min_{d \in \mathbb{R}^p, \alpha \in \mathbb{R}} \alpha + \frac{1}{2} \|d\|_2^2 \\
        \text{s.t.} & \quad \nabla \mathcal{L}_{\text{vol}}(\theta_t)^T d \leq \alpha, \\
        & \quad \nabla \mathcal{L}_{\text{feas}}(\theta_t)^T d \leq \alpha.
    \end{aligned}
    \end{equation}

    \item If $\mathcal{L}_{\text{vol}}(\theta_t) < \beta \mathcal{L}_{\text{feas}}(\theta_t) + \epsilon_{\text{lb}}$:
    \begin{equation}
    \begin{aligned}
        \label{eq:primal_3}
        (d_t, \alpha_t) = \arg & \min_{d \in \mathbb{R}^p, \alpha \in \mathbb{R}} \alpha + \frac{1}{2} \|d\|_2^2 \\
        \text{s.t.} & \quad [\nabla (\beta \mathcal{L}_{\text{feas}}(\theta_t) - \mathcal{L}_{\text{vol}}(\theta_t))]^T d \leq \alpha, \\
        & \quad \nabla \mathcal{L}_{\text{feas}}(\theta_t)^T d \leq \alpha.
    \end{aligned}
    \end{equation}
\end{enumerate}

\begin{figure}[ht]
    \centering
    \includegraphics[width=.7\linewidth]{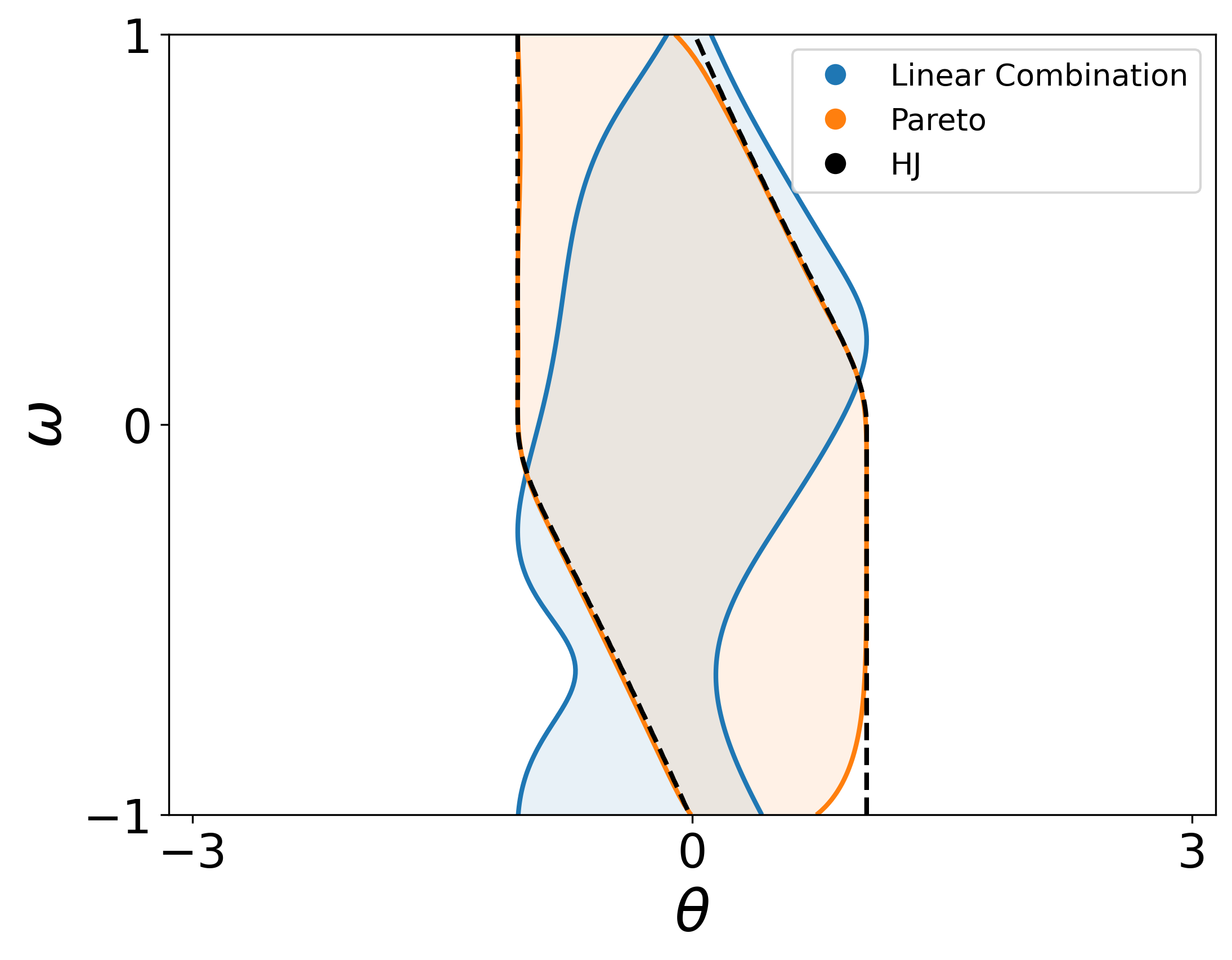}
    \caption{The inner safe sets of the inverted pendulum problem obtained by the LCCBF (blue solid line) and PCBF (orange solid line) methods, alongside the infinite-time viability kernel (black dashed line) computed using HJ reachability. The region enclosed by the black dashed line represents the true safe states.}
    \label{fig:IP_full}
\end{figure}

\begin{table}[htbp]
    \centering
    \begin{tabular}{lccc}
        \toprule
        Model & Sampling Method & Sampling Num & Training Time (s) \\
        \midrule
        LCCBF & Boundary & 10000 & 17823 $\pm$ 713 \\
        LCCBF & Gaussian & 10000 & \textbf{1143 $\pm$ 58} \\
        PCBF & Boundary & 10000 & 18356 $\pm$ 847 \\
        PCBF & Gaussian & 10000 & \textbf{1207 $\pm$ 69} \\
        \bottomrule
    \end{tabular}
    \caption{Training Time for Different Sampling Methods on NVIDIA GeForce RTX 3080ti}
    \label{tab:training_times}
\end{table}

\begin{algorithm*}[ht]
\caption{Pareto Control Barrier Function (PCBF)}
\label{pseudocode}
\begin{algorithmic}[1]

\Function{SolveBaseSubproblem}{$\mathcal{L}_1, \mathcal{L}_2$}
\Comment{See Appendix \ref{appendix:analytical-solution}}

    $\lambda \gets \max\left(\min\left( 
    \frac{(\nabla \mathcal{L}_2 - \nabla \mathcal{L}_1)^T \nabla \mathcal{L}_2}
    {\|\nabla \mathcal{L}_1 - \nabla \mathcal{L}_2\|_2^2}, 1\right), 0\right)$

    \State \Return $-\lambda \nabla \mathcal{L}_1 - (1-\lambda) \nabla \mathcal{L}_2$
\EndFunction

\Require Initial parameter $\theta_0 \in \mathbb{R}^p$, learning rate $\eta > 0$, maximum iterations $T \in \mathbb{N}$, hyperparameters $\beta, \epsilon_{\text{ub}}, \epsilon_{\text{lb}}$, regularization weight $\gamma > 0$, safe state $x_e$, dataset $\mathbb{X}_{\mathcal{N}}$.

\For{$t = 0, 1, \ldots, T-1$}
    \State
    $\mathcal{L}_f \gets \mathcal{L}_{\text{feas}}(\theta_t)$, 
    $\mathcal{L}_v \gets \mathcal{L}_{\text{vol}}(\theta_t)$
    
    \If{$\mathcal{L}_v > \beta\mathcal{L}_f + \epsilon_{\text{ub}}$}
        \State $d_t \gets \Call{SolveBaseSubproblem}{\mathcal{L}_v - \beta \mathcal{L}_f, \mathcal{L}_v}$
        \Comment{Eq. \eqref{subproblem-up-dt}, \eqref{subproblem-up-lambda}}
        
    \ElsIf{$\beta\mathcal{L}_f + \epsilon_{\text{lb}} \leq \mathcal{L}_v \leq \beta\mathcal{L}_f + \epsilon_{\text{ub}}$}
        \State $d_t \gets \Call{SolveBaseSubproblem}{\mathcal{L}_f, \mathcal{L}_v}$
        \Comment{Eq. \eqref{subproblem-middle-dt}, \eqref{subproblem-middle-lambda}}
        
    \Else
        \State $d_t \gets \Call{SolveBaseSubproblem}{\beta \mathcal{L}_f - \mathcal{L}_v, \mathcal{L}_f}$
        \Comment{Eq. \eqref{subproblem-below-dt}, \eqref{subproblem-below-lambda}}
    \EndIf
    
    \State 
    $g_t \gets d_t - \gamma\nabla\mathcal{L}_{\text{vol}}(\theta_t)$
    \Comment{Apply regularization}
    \State
    $\theta_{t+1} \gets \theta_t + \eta g_t$
\EndFor
\State \Return $\theta_T$
\end{algorithmic}
\end{algorithm*}

All three subproblems are instances of problem \eqref{eq:problem-base} with \( M=2 \). According to lemma \ref{lemma:steepest-descent}, both losses will decrease as \( \theta_t \) approaches a Pareto critical solution. 
For subproblem \eqref{eq:primal_1}, the update direction \( d_t \) reduces both \( \mathcal{L}_{\text{feas}} \) and \( \mathcal{L}_{\text{vol}} - \beta \mathcal{L}_{\text{feas}} \), where the second term measures the violation of the upper bound. Once this violation decreases to be less than or equal to \( \epsilon_{\text{ub}} \), the loss vector \( \mathbf{L}(\theta_t) \) returns to within the bounds.
Similarly, subproblem \eqref{eq:primal_3} prevents \( \theta_t \) from violating the lower bound.
For subproblem \eqref{eq:primal_2}, the constraints incorporate both feasibility and volume losses. Following lemma \ref{lemma:steepest-descent}, both losses decrease until a Pareto critical solution is reached.

\begin{figure*}[ht]
    \centering
    \includegraphics[width=\linewidth]{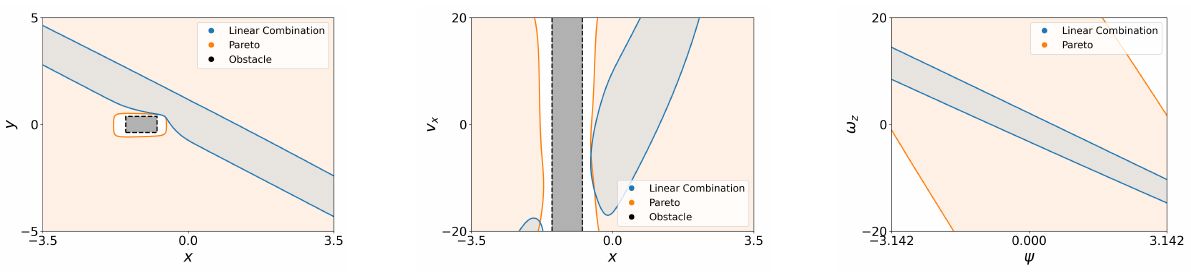}
    \caption{The left, middle, and right subfigures present 2D slices of the learned inner safe sets obtained by the PCBF and LCCBF methods corresponding to the $x-y$, $x-v_x$, and $\psi-\omega_z$ planes. All other states are set to zero in each slice. In the left and middle subfigures, the obstacle's projection in the corresponding 2D plane is indicated by black dashed lines, representing regions unsafe for the quadrotor.}
    \label{fig:combined_figures}
\end{figure*}

To avoid \( \theta_t \) converging to a Pareto critical solution that is not Pareto optimal, a regularization term \( \gamma \mathcal{L}_{\text{vol}} \) is introduced, where \( \gamma \) is a small positive constant. This term serves two purposes: preventing convergence to a suboptimal Pareto critical solution and accelerating the training process. 
The pseudocode of the PCBF algorithm is shown in Algorithm \ref{pseudocode}, where analytical solutions for all subproblems are provided; for the proof, please refer to the Appendix.

\section{SIMULATIONS}
\label{sec:experiments}

In this section, the performance of the PCBF algorithm is evaluated. First, a 2D inverted pendulum example is employed, where the low dimensionality allows HJ Reachability to compute the infinite-time viability kernel as the ground truth for the largest inner safe set. Comparison of the PCBF results with HJ validates its effectiveness. Next, simulations were conducted in a 12-dimensional quadrotor obstacle avoidance scenario to assess PCBF's performance in this more complex environment. Neural CBFs trained with linear combination loss (LCCBF) serve as our baseline. All neural networks are trained until convergence and share the same architecture (3 layers with 256 neurons and tanh activations).

\subsection{Inverted Pendulum}

The inverted pendulum dynamics are $\dot{\theta} = \omega, \dot{\omega} = \sin(\theta) + u$, where $\theta$ is the angle of the pendulum from the upright position, $\omega$ is the angular velocity and $u$ is the control input. The safe set is defined as $|\theta| \leq \frac{\pi}{3}$, and the input constraint is $|u| \leq 1$.

Figure \ref{fig:IP_full} shows a strong overlap between the inner safe set generated by PCBF and the infinite-time viability kernel from HJ reachability. Almost all states in the PCBF inner safe set fall within the viability kernel, demonstrating that PCBF effectively maximizes the safe set while ensuring safety. In contrast, the LCCBF-generated set shows significant discrepancies, with many states outside the kernel, leading to potential safety violations in practice.

The training times of different methods are reported in Table \ref{tab:training_times}. In the sampling method column, "Boundary" refers to the sampling approach used in \cite{liu2023safe}, while "Gaussian" represents the method proposed in this study. The results demonstrate a significant improvement in training efficiency with our new approach.

\subsection{Quadrotor}

The state of the quadrotor involves its position $x, y, z$, velocity $ v_x, v_y, v_z$, Euler angles $\phi, \theta, \psi$, and angular velocity $\omega_x, \omega_y, \omega_z$. The state dimension is 12 and the control dimension is 4. The system dynamics and input constraints remain the same as in \cite{beard2008quadrotor}. A rectangular obstacle is placed in the scenario and the safe set is defined such that the quadrotor must maintain a minimum distance from the obstacle. Specifically, the safe set is \(x \leq -1.5\) or \(x \geq -0.75\) or \(y \leq -0.375\) or \(y \geq 0.375\).

Figure \ref{fig:combined_figures} illustrates the inner safe sets of PCBF and LCCBF in various dimensions through 2D slices. Notably, the PCBF consistently yields a significantly larger inner safe set slice compared to the LCCBF across all slices. For example, in the left subfigure in Figure \ref{fig:combined_figures}, where all states except $x$ and $y$ are set to zero, the quadrotor can hover in any unobstructed region. Consequently, the area outside the black region is safe. The inner safe set of the PCBF encompasses a substantial portion of the safe state, while the LCCBF's inner safe set is confined to a relatively small area. The volume of the inner safe set for the LCCBF is only 14\% of that for the PCBF in the entire 12-dimensional state space.

\begin{figure}[ht]
    \centering
    \includegraphics[width=.8\linewidth]{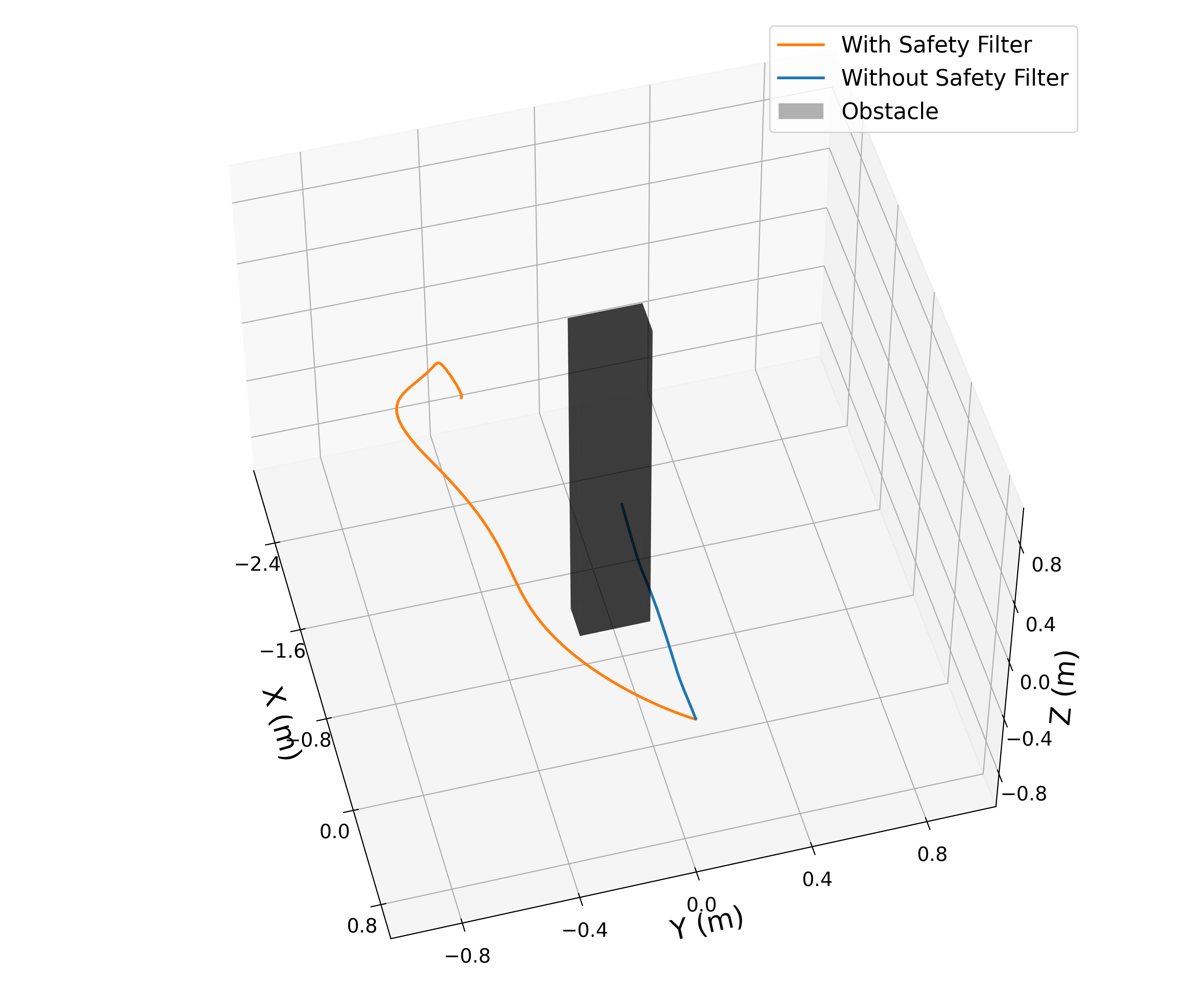}
    \caption{Trajectories of the quadrotor with and without the PCBF safety filter in the presence of a rectangle obstacle.}
    \label{fig:trajectory}
\end{figure}

Simulations are conducted using RotorPy \cite{folk2023rotorpy} in a windless environment with a refresh rate of 300 Hz. The nominal policy is based on the control algorithm proposed by \cite{lee2010geometric}.  When the PCBF safety filter is activated, problem \ref{eq:cbf-qp} is solved to derive a safe control policy.
Figure \ref{fig:trajectory} illustrates the simulation results, highlighting the quadrotor trajectories with and without the PCBF safety filter. The trajectory without the safety filter collides with an obstacle, while the trajectory using the safety filter successfully avoids it. 

\section{CONCLUSIONS}
\label{sec:conclusions}

In this paper, the Pareto Control Barrier Function (PCBF) algorithm is introduced to maximize the inner safe set under input constraints. The PCBF algorithm significantly outperforms traditional LCCBF methods. Comparisons with HJ reachability for the inverted pendulum showed that PCBF closely approximates the largest inner safe set. Simulations with the quadrotor system further demonstrated PCBF's effectiveness in complex, high-dimensional environments. Future work will focus on adapting the PCBF framework to handle uncertainties in system dynamics and external disturbances.








\bibliography{IEEEabrv,ref}
\bibliographystyle{IEEEtran}

\clearpage

\section*{APPENDIX}
\subsection{Proof of lemma \ref{lemma:ism-pareto}}

Let $\theta^*$ be the solution to problem \ref{eq:problem-ism-loss}. Suppose $\theta^*$ is not Pareto optimal. Then there exists a Pareto optimal $\hat{\theta}$ such that:

\begin{align}
    \mathcal{L}_{\text{feas}}(\hat{\theta}) &\leq \mathcal{L}_{\text{feas}}(\theta^*) = 0 \\
    \mathcal{L}_{\text{vol}}(\hat{\theta}) &\leq \mathcal{L}_{\text{vol}}(\theta^*)
\end{align}

\noindent with at least one of the inequalities being strict. However, since $\mathcal{L}_{\text{feas}}(\theta^*) = 0$, we must have $\mathcal{L}_{\text{feas}}(\hat{\theta}) = 0$ as well. This means $\hat{\theta}$ is also a feasible solution to problem \ref{eq:problem-ism-loss}.
If $\mathcal{L}_{\text{vol}}(\hat{\theta}) < \mathcal{L}_{\text{vol}}(\theta^*)$, then $\hat{\theta}$ would be a better solution to problem \ref{eq:problem-ism-loss} than $\theta^*$, which contradicts the assumption that $\theta^*$ is the solution to this problem.
If $\mathcal{L}_{\text{vol}}(\hat{\theta}) = \mathcal{L}_{\text{vol}}(\theta^*)$, then $\hat{\theta}$ is also an optimal solution to problem \ref{eq:problem-ism-loss}, and thus satisfies the conditions of the lemma.

Therefore, our assumption must be false, and $\theta^*$ is indeed Pareto optimal.
By definition, $[\mathcal{L}_{\text{feas}}(\theta^*), \mathcal{L}_{\text{vol}}(\theta^*)]^T$ must lie on the Pareto front.

\subsection{Analytical solution for PCBF subproblems}
\label{appendix:analytical-solution}
Subproblems \eqref{eq:primal_1}, \eqref{eq:primal_2}, and \eqref{eq:primal_3} are special cases of problem \eqref{eq:problem-base} with $M = 2$. We first provide the general solution for problem \eqref{eq:problem-base}. For $M = 2$, the problem becomes:

\begin{equation}
\begin{aligned}
(d_t, \alpha_t) = \arg & \min_{d \in \mathbb{R}^p, \alpha \in \mathbb{R}} \alpha + \frac{1}{2}\|d\|^2 \\
\text{s.t.} & \quad \nabla \mathcal{L}_1(\theta_t)^T d \leq \alpha \\
            & \quad \nabla \mathcal{L}_2(\theta_t)^T d \leq \alpha
\end{aligned}
\label{eq:problem-base-M2}
\end{equation}

As problem \eqref{eq:problem-base-M2} is a QP problem, we can apply the Karush–Kuhn–Tucker (KKT) conditions:
    \begin{align}
    d_t + \lambda_1 \nabla \mathcal{L}_1(\theta_t) + \lambda_2 \nabla \mathcal{L}_2(\theta_t) = 0 \\
    1 - \lambda_1 - \lambda_2 = 0 \\
    \lambda_1(\nabla \mathcal{L}_1(\theta_t)^T d_t - \alpha_t) = 0 \\
    \lambda_2(\nabla \mathcal{L}_2(\theta_t)^T d_t - \alpha_t) = 0 \\
    \lambda_1, \lambda_2 \geq 0 \\
    \nabla \mathcal{L}_1(\theta_t)^T d_t - \alpha_t \leq 0 \\
    \nabla \mathcal{L}_2(\theta_t)^T d_t - \alpha_t \leq 0
    \end{align}

    From the first two conditions, we get:
    \begin{align}
        d_t = -\lambda_1 \nabla \mathcal{L}_1(\theta_t) - \lambda_2 \nabla \mathcal{L}_2(\theta_t) \\
        \lambda_1 + \lambda_2 = 1
    \end{align}

    Let $\lambda_1 = \lambda$ and $\lambda_2 = 1 - \lambda$. Then:
    \begin{equation}
        d_t = -\lambda \nabla \mathcal{L}_1(\theta_t) - (1-\lambda) \nabla \mathcal{L}_2(\theta_t)
    \end{equation}

    To find $\lambda$, we use the complementary slackness conditions:
    \begin{align}
        \lambda(\nabla \mathcal{L}_1(\theta_t)^T d_t - \alpha_t) &= 0 \\
        (1-\lambda)(\nabla \mathcal{L}_2(\theta_t)^T d_t - \alpha_t) &= 0
    \end{align}

    This leads to three cases:
    \begin{itemize}
        \item If $0 < \lambda < 1$, then $\nabla \mathcal{L}_1(\theta_t)^T d_t = \nabla \mathcal{L}_2(\theta_t)^T d_t = \alpha_t$
        Substituting $d_t$ and solving for $\lambda$:

        \begin{equation}
        \begin{aligned}
        &-\lambda \|\nabla \mathcal{L}_1(\theta_t)\|^2 - (1-\lambda) \nabla \mathcal{L}_1(\theta_t)^T \nabla \mathcal{L}_2(\theta_t) \\
        &= -\lambda \nabla \mathcal{L}_2(\theta_t)^T \nabla \mathcal{L}_1(\theta_t) - (1-\lambda) \|\nabla \mathcal{L}_2(\theta_t)\|^2
        \end{aligned}
        \end{equation}
        
        \begin{equation}
        \lambda = \frac{(\nabla \mathcal{L}_2(\theta_t) - \nabla \mathcal{L}_1(\theta_t))^T \nabla \mathcal{L}_2(\theta_t)}
        {\|\nabla \mathcal{L}_1(\theta_t) - \nabla \mathcal{L}_2(\theta_t)\|_2^2}
        \end{equation}
        
        \item If $\lambda = 0$, then $\nabla \mathcal{L}_2(\theta_t)^T d_t = \alpha_t$ and $\nabla \mathcal{L}_1(\theta_t)^T d_t \leq \alpha_t$
        \item If $\lambda = 1$, then $\nabla \mathcal{L}_1(\theta_t)^T d = \alpha_t$ and $\nabla \mathcal{L}_2(\theta_t)^T d \leq \alpha_t$
    \end{itemize}

    Combining these conditions leads to the following solution:
    \begin{equation}
        \lambda = \max\left(\min\left( 
        \frac{(\nabla \mathcal{L}_2(\theta_t) - \nabla \mathcal{L}_1(\theta_t))^T \nabla \mathcal{L}_2(\theta_t)}
        {\|\nabla \mathcal{L}_1(\theta_t) - \nabla \mathcal{L}_2(\theta_t)\|_2^2}, 1\right), 0\right)
    \end{equation}

Replacing $\mathcal{L}_1$ and $\mathcal{L}_2$ with specific losses in problems \eqref{eq:primal_1}, \eqref{eq:primal_2}, and \eqref{eq:primal_3} leads to the following analytical solutions:

\begin{enumerate}
    \item If $\mathcal{L}_{\text{vol}}(\theta_t) > \beta \mathcal{L}_{\text{feas}}(\theta_t) + \epsilon_{\text{ub}}$:
    \begin{equation}
        d_t = (\lambda + \lambda\beta - 1)\nabla \mathcal{L}_\text{feas}(\theta_t) - \lambda\nabla \mathcal{L}_\text{vol}(\theta_t)
    \label{subproblem-up-dt}
    \end{equation}
    where
    \begin{equation}
    \begin{aligned}
        \lambda = \max\Bigg(\min\Bigg( &
        \frac{((1 + \beta)\nabla \mathcal{L}_\text{feas}(\theta_t) - \nabla \mathcal{L}_\text{vol}(\theta_t))^T}
        {\|(1 + \beta)\nabla \mathcal{L}_\text{feas}(\theta_t) - \nabla \mathcal{L}_\text{vol}(\theta_t)\|_2^2} \\
        & \cdot \nabla \mathcal{L}_\text{feas}(\theta_t), \,
        1\Bigg), 0\Bigg)
    \label{subproblem-up-lambda}
    \end{aligned}
    \end{equation}
        
    \item If $\beta \mathcal{L}_{\text{feas}}(\theta_t) + \epsilon_{\text{lb}} \leq \mathcal{L}_{\text{vol}}(\theta_t) \leq \beta \mathcal{L}_{\text{feas}}(\theta_t) + \epsilon_{\text{ub}}$:
    \begin{equation}
        d_t = -\lambda\nabla \mathcal{L}_\text{vol}(\theta_t) - (1 - \lambda)\nabla \mathcal{L}_\text{feas}(\theta_t)
    \label{subproblem-middle-dt}
    \end{equation}
    where
    \begin{equation}
    \begin{aligned}
        \lambda = \max\Bigg(\min\Bigg( &
        \frac{(\nabla \mathcal{L}_\text{feas}(\theta_t) - \nabla \mathcal{L}_\text{vol}(\theta_t))^T \nabla \mathcal{L}_\text{feas}(\theta_t)}
        {\|\nabla \mathcal{L}_\text{vol}(\theta_t) - \nabla \mathcal{L}_\text{feas}(\theta_t)\|_2^2}, \\
        & 1\Bigg), 0\Bigg)
    \label{subproblem-middle-lambda}
    \end{aligned}
    \end{equation}
    
    \item If $\mathcal{L}_{\text{vol}}(\theta_t) < \beta \mathcal{L}_{\text{feas}}(\theta_t) + \epsilon_{\text{lb}}$:
    \begin{equation}
        d_t = (\lambda - \lambda\beta - 1)\nabla \mathcal{L}_\text{feas}(\theta_t) + \lambda\nabla \mathcal{L}_\text{vol}(\theta_t)
    \label{subproblem-below-dt}
    \end{equation}
    where
    \begin{equation}
    \begin{aligned}
        \lambda = \max\Bigg(\min\Bigg( &
        \frac{(\nabla \mathcal{L}_\text{vol}(\theta_t) + (1 - \beta)\nabla \mathcal{L}_\text{feas}(\theta_t))^T}
        {\|\nabla \mathcal{L}_\text{vol}(\theta_t) + (1 - \beta)\nabla \mathcal{L}_\text{feas}(\theta_t)\|_2^2} \\
        & \cdot \nabla \mathcal{L}_\text{feas}(\theta_t), \,
        1\Bigg), 0\Bigg)
    \label{subproblem-below-lambda}
    \end{aligned}
    \end{equation}
\end{enumerate}

\end{document}